\documentclass[draft]{amsart}
\usepackage{amscd}

\begin{document}
\title [A note on finite Blaschke products ] {A note on the characterization of finite Blaschke products}
\author{Nina Zorboska *}

\address{Department of Mathematics, University of Manitoba, Winnipeg, MB, R3T 2N2} 
\email{zorbosk@cc.umanitoba.ca}
\thanks{* Research supported in part by NSERC grant.}
\date{December 2015}

\begin{abstract}
We give a slight generalization of the characterization of finite Blaschke products given in \cite{Zorb1}. The characterization uses the boundary behaviour of a weighted local hyperbolic distortion of an analytic self-map of the unit disk.
\end{abstract}

\maketitle

\bigskip

\def\C{\mathbb{C}}
\def\D{\mathbb{D}}
\def\a{\alpha}
\def\b{\beta}
\def\g{\gamma}
\def\d{\delta}
\def\e{\varepsilon}
\def\l{\lambda}
\def\r{\rho}
\def\s{\sigma}
\def\t{\theta}
\def\B{\mathcal{B}}
\def\p{\phi}
\def\vp{\varphi}


\bigskip

The main purpose of this note is to provide a slight generalization of the characterization of finite Blaschke products given in \cite{Zorb1}.
The new result includes a wider range of weights than the original one. The notation and the main methods follow \cite{Zorb1}. For completeness, we include a short list of definitions and basic concepts used in the statement and in the proof of the result. 

Let $\D$ be the open unit disk in the complex plane $\C$. 
For a non-constant analytic function $\phi$ that maps the unit disk into itself and for $\alpha>0$, let
$$\tau_{\p, \a}(z)=\frac{(1-|z|^2)^{\a}|\p'(z)|}{(1-|\p(z)|^2)^\a}.$$
We will say that $\tau_{\p, \a}(z)$ is the local \textit{$\alpha$-hyperbolic distortion} of $\p$ at $z$. 

The motivation for this definition comes from the classical hyperbolic case with $\alpha=1$.
Recall that for $z \in \D$, $\l(z)=\frac{1}{1-|z|^2}$ is the density of the hyperbolic metric on $\D$, 
and that for an analytic, non-constant map
$\p: \D \to \D$, the pull-back of the hyperbolic metric is defined by $\p(\l)^*(z) = \frac{|\p' (z)}{1-|\p(z)|^2}$. Thus, 
$$\tau_\p(z)= \tau_{\p, 1}(z)=\frac{(1-|z|^2)|\p'(z)|}{(1-|\p(z)|^2)} = \frac{\p(\l)^*(z)}{\l(z)}$$
is the usual local hyperbolic distortion of $\p$ at $z$. 
Similarly, for $\alpha>0$, we can think of $\tau_{\p, \a}(z)$, the local $\alpha$-hyperbolic distortion of $\p$ at $z$, as the pull-back
by $\p$ of the $\alpha$-hyperbolic metric on $\D$ with density $\l_{\a}(z)=\frac{1}{(1-|z|^2)^{\a}}$.
(See, for example,  \cite{BeMi} for further details on these and other related basic notions, results and references).

\medskip

By the classical Schwarz-Pick lemma, $\tau_\p(z) \leq 1$ for all $z\in \D$, i.e. every self-map of the unit disk is a hyperbolic contraction.
Furthermore, if the equality holds for one $z \in \D$, then $\p$ is a disk automorphism, and so the equality must hold for every $z \in \D$. 
Thus, the maximal possible hyperbolic distortion is attained in the disk only when the map is a disk automorphism. The situation is very  different when $\a \ne 1$. For example, there are analytic self-maps of $\D$ for which $\tau_{\p, \a}(z)$ is not even bounded when 
$0<\a<1$. For more details and further references on this, see \cite{Zorb1}.

\medskip

In 1986, M. Heins used in \cite{He} the boundary behaviour of the hyperbolic distortion $\tau_\p$ to characterize the finite Blaschke products among the class of analytic self-maps of the unit disk.

 \medskip

  \textbf{Theorem A. (\cite{He}):} Let $\p$ be an analytic self-map of $\D$. Then $\p$ is a finite Blaschke product if and only if
  $\lim_{|z| \to 1} \tau_\p (z) = 1.$
 
 \medskip

 A more recent result of  D. Kraus, O. Roth and  S. Rucheweyh  from 2007 generalized Heins' result to a characterization of analytic
 boundary behaviour of self-maps of $\D$ on subarcs of the unit circle. 
  
 \medskip
 
 \textbf{Theorem B. (\cite{KRR}):} Let $\p$ be an analytic self-map of $\D$ and let $\Gamma$ be an open subarc of $\partial \D$.
 Then the following are equivalent:
 
  (a) For every $\zeta \in \Gamma$, $\liminf_{z \to \zeta} \tau_\p (z) >0$.
  
  (b) For every $\zeta \in \Gamma$, $\lim_{z \to \zeta} \tau_\p (z) =1$.
  
  (c) $\p$ has an analytic extension across $\Gamma$ and $\p(\Gamma) \subset \partial\D$.
 
 \medskip

They also asked questions about possible generalizations of these results in terms of angular, i.e. non-tangential, limits. 
Recall that for $\zeta \in \partial \D$ and $\gamma>1$, the angular (or non-tangential) region $\Gamma_{\g}(\zeta)$ in $\D$
is defined by
$$\Gamma_{\g}(\zeta) = \{z \in \D: |\zeta - z|  \le \g(1-|z|^2)\}.$$ 
If $z \to \zeta$ through the angular region $\Gamma_{\g}(\zeta)$, we say that the corresponding limit is an 
{\it{angular limit}}. This type of limit will be denoted by $\angle \lim_{z \to \zeta}$. Recall that since an analytic self-map $\p$ of the unit disk is in $H^{\infty}(\D)$, it has radial (and angular) limits almost everywhere on the unit circle. We will denote the radial extension function with the same symbol $\p$.

Moreover, an analytic map $\p: \D \to \D$ has an 
$\it {angular ~derivative}$ at $\zeta \in \partial \D$ if there exists $\xi \in \partial \D$ such that the angular limit
$$\angle \lim_{z \to \zeta} \frac{\p(z) - \xi}{z - \zeta}$$ exists. In this case the value of the limit 
is called an angular derivative of $\p$ at $\zeta$, and will be denoted by $\p'(\zeta)$. 
By the Julia-Carath\'eodory theorem (see \cite{Sh}), the existence of the angular derivative at $\zeta$ is equivalent to 
$$0< \liminf_{z \to \zeta}\frac{1-|\p(z)|}{1-|z|}=|\p'(\zeta)| < \infty,$$
with the limit infimum attained in an angular approach to $\zeta$.

\medskip

In 2013, D. Kraus  specified in \cite{Kr} the class of functions characterized by the existence of (nonzero) angular boundary limits 
of their hyperbolic distortion a.e. on the unit circle. 
 
 \medskip
 
 \textbf{Theorem C.(\cite{Kr}):} If $\p$ is an analytic self-map of $\D$, then $\angle \lim_{z \to \zeta} \tau_\p (z) =1$ 
 for almost every $\zeta \in \partial\D$ if and only if $\p$ is an inner function with finite angular derivatives 
 at almost every point in $\partial\D$.  
  
 \medskip

 Since there exist infinite Blaschke products
with finite angular derivatives a.e. on $\partial\D$, the class of functions determined by the last theorem is much larger than the class of finite Blaschke products. Hence, the angular limits boundary condition on the classical hyperbolic distortion $\tau_\p$ will not determine the class of finite Blaschke products. We were able to show in \cite{Zorb1} that we can accomplish this if, instead, we use angular limit boundary conditions on $\tau_{\p, \a}$,  with $ \a>1$. 
The following theorem shows that actually, this result can be generalized to any $\a \ne 1$. Parts of the proof of the theorem are the same as in \cite{Zorb1}. We provide some of them here for completeness.

 \medskip

\textbf{Theorem 1. } Let $\p$ be a non-constant self-map of $\D$. Then $\p$ is a
finite Blaschke product if and only if there exist $\a \ne 1$ and $c>0$ such that $\tau_{\p, \a}$ is bounded
and such that for almost every $\zeta \in \partial\D$
$$\angle \liminf_{z \to \zeta} \tau_{\p, \a}(z) \geq c.$$
\begin{proof}
If $\p$ is a finite Blaschke product, then $\tau_{\p, \a}$ is bounded for any $\a>0$ as shown, for example, in \cite{Zorb1}. Since $\p$ has an analytic extension across the unit circle, and
$|\p(\zeta)|=1$ for all $\zeta \in \partial\D$, by the Julia-Carath\'eodory theorem  $\p$ has an angular derivative equal to the regular derivative of $\p$ at $\zeta$.
Using again the Julia-Carath\'eodory theorem characterization of the angular derivative, one gets that 
$$\angle \lim_{z \to \zeta}  \tau_{\p, \a}(z) = |\p' (\zeta)|^{1-\a}.$$  
If $\a>1$, since $\p'$ is continuous on $\overline{\D}$, we get that 
$$\angle \lim_{z \to \zeta}  \tau_{\p, \a}(z) = |\p' (\zeta)|^{1-\a} \geq || \p' ||_{\infty}^{1-\a} >0.$$
If $0<\a<1$, since by the Schwartz Lemma $\frac{1-|\p(z)|^2}{1-|z|^2} \geq \frac{1-|\p(0)|}{1+|\p(0)|}$, we get that
$$\angle \lim_{z \to \zeta}  \tau_{\p, \a}(z) = |\p' (\zeta)|^{1-\a}= \angle \lim_{z \to \zeta} \left( \frac{1-|\p(z)|^2}{1-|z|^2} \right)^{1-\a}
\ge \left(\frac{1-|\p(0)|}{1+|\p(0)|}\right)^{1-\a}>0.$$

\medskip

For the other direction of the proof, we will show first, by using a similar idea as in  the proof of Theorem 2.3 from \cite{Kr}, that if 
$\exists \a>0, c>0$ such that 
$\angle \liminf_{z \to \zeta} \tau_{\p, \a}(z) \geq c$ for almost every $\zeta \in \D$, then $\p$ has to be an inner function.
Note that $1-|z|^2 \to 0$ as $z \to \zeta \in \partial \D$, and so for almost all $\zeta$ in $\partial \D$
$$\angle \lim_{z \to \zeta} \frac{|\p'(z)|}{(1-|\p(z)|)^\a} = \infty.$$
Thus, if $\angle \lim_{z \to \zeta} |\p(z)|=r<1$, then 
$\angle \lim_{z \to \zeta} |\p'(z)|= \infty$. But by Privalov's Theorem \cite[p.~140]{Po}, this can happen only for $\zeta$ in a set of measure zero, i.e.
$\angle \lim_{z \to \zeta} |\p(z)|=1$ for almost every $\zeta \in \partial\D$, and so $\p$ is inner.

In the case when $0<\a<1$, since $\p$ is inner and $\tau_{\p, \a}$ is bounded, $\p$ must be a finite Blaschke product. This follows,
for example, from the fact that then $\p$ has to be in the disk algebra (see \cite[p.~74]{Du}). 

When $\a>1$, and for almost every $\zeta \in \partial\D$
$$\angle \liminf_{z \to \zeta} \tau_{\p, \a}(z) = \angle \liminf_{z \to \zeta} \left( \frac{1-|z|^2}{1-|\p(z)|^2}\right) ^{\a-1}\tau_{\p}(z) \geq c >0,$$
then, since by the Schwarz-Pick Lemma $\tau_{\p}(z) \leq 1$ for every $z \in \D$, and since $\a-1 >0$, we have that 
$$\angle \liminf_{z \to \zeta}\frac{1-|z|^2}{1-|\p(z)|^2} \geq c^ \frac{1}{\a-1}.$$
Thus, for almost every $\zeta$ in $\partial\D$
$$\angle \limsup_{z \to \zeta}\frac{1-|\p(z)|^2}{1-|z|^2} \leq  \left(\frac{1}{c}\right)^{\frac {1}{\a-1}},$$
and so, for almost every $\zeta$ in $\partial\D$ it must be that
$$|\p'(\zeta)| =  \liminf_{z \to \zeta}\frac{1-|\p(z)|}{1-|z|} \leq 2 \left(\frac{1}{c}\right)^{\frac {1}{\a-1}} < \infty.$$ 
Hence,  we have that on one hand $\p$ is inner, as shown above, and on the other hand for almost every $\zeta$ in $\partial \D$ the modulus of the angular derivative at $\zeta$ is bounded above by a constant depending only on $\a$ and $c$.

By a result of P. Ahern and D. Clark (see \cite{AhCl}), for an inner function $\p$ we have that $|\p'|$ is in 
$L^p(\partial\D), 0< p \leq \infty$ if and only if $\p'$ is in $H^p(\D)$. Furthermore, if
$\p' \in H^{1}(\D)$, then $\p$ must be a finite Blaschke product.
Since in our case  $|\p'|$ is almost everywhere bounded by $2 \left(\frac{1}{c}\right)^{\frac {1}{\a-1}}$ on $\partial\D$, 
by \cite{AhCl} we have that $\p'$ belongs to $H^{\infty}(\D)$, and so $\p$ must be a finite Blaschke product.
\end{proof}

\medskip

\textbf{Remark}: Note that as a consequence of the theorem, one also gets an interesting result reflecting on the rigidity of the local hyperbolic distortion. Namely, if there exist $\a \ne 1$ and $c>0$ such that $\tau_{\p, \a}$ is bounded and such that for almost 
every $\zeta \in \partial\D$ we have that
$\angle \liminf_{z \to \zeta} \tau_{\p, \a}(z) \geq c,$
then  $\angle \lim_{z \to \zeta} \tau_{\p, \a}(z)$ exists for every $\zeta$ in $\partial\D$ and equals $|\p'(\zeta)|^{1-\a}$.
Also, the proof of the theorem shows that if in the previous statement we replace $\partial\D$ by a subset $E$ of positive Lebesque (linear) measure, then $\angle \lim_{z \to \zeta} |\p(z)|=1$ for almost every $\zeta \in E$ and the conclusion on the rigidity of  
$\tau_{\p, \a}$ holds almost everywhere on $E$.

\bigskip

\bigskip

\end{document}